\documentclass[phd,tocprelim]{cornell}

\let\ifpdf\relax
\usepackage{latexsym,amsfonts,amsmath,amssymb,mathabx,enumitem,graphicx,fancyhdr,amsthm,pstricks,graphics,moreverb,subfigure,epsfig,subfigure,hangcaption,txfonts,palatino,lipsum,mathabx}
\usepackage[english]{babel}
\usepackage[autostyle]{csquotes}
\usepackage[all]{xy}
\usepackage{mathtools} 

\makeatletter
\newcommand{\leqnomode}{\tagsleft@true}
\newcommand{\reqnomode}{\tagsleft@false}
\makeatother
\renewcommand{\labelenumii}

\newcommand{\Z}{{\mathbb{Z}}}

\newtheorem{theorem}{Theorem}
\newtheorem{corollary}{Corollary}
\newtheorem*{claim*}{Claim}
\theoremstyle{definition}
\newtheorem{definition}{Definition}

\theoremstyle{remark}

\DeclareMathOperator{\rank}{rank}

\renewcommand{\caption}[1]{\singlespacing\hangcaption{#1}\normalspacing}

\title {A Study of finitely generated Free Groups via the Fundamental Groups}
\author {\textbf{Gongping Niu}}
\conferraldate {June}{2017}
\degreefield {}
\copyrightholder{Gongping Niu}
\copyrightyear{2017}

\begin{document}
	
	\maketitle
	\makecopyright
	
	\begin{abstract}
		Free groups have many applications in Algebraic Topology. In this paper I specifically study the finitely generated free groups by using the covering spaces and fundamental groups. By the Van Kampen's theorem, we have a famous fact that $\pi_1(S^1\vee S^1)\backsimeq \Z * \Z=F_2$. Therefore, to study $F_n$, we could try to figure out the covering spaces of $S^1\vee S^1$ or even $\bigvee_i S^1$. And in the appendix B, we  prove the Nielsen-Schreier theorem which we will use this to study finitely index subgroups of $F_n$.
	\end{abstract}
	
	\contentspage
	
	\normalspacing \setcounter{page}{1} \pagenumbering{arabic}
	\pagestyle{cornell} \addtolength{\parskip}{0.5\baselineskip}
		\reqnomode
	Free groups have many applications in Algebraic Topology. In this paper I specifically study the finitely generated free groups by using the covering spaces and fundamental groups. To review the definition of Free groups, you can read \cite{Dummit} or the Appendix A. 
	
	By the Van Kampen's theorem, we have a famous fact that $\pi_1(S^1\vee S^1)\backsimeq \Z * \Z=F_2$. Therefore, to study $F_n$, we could try to figure out the covering spaces of $S^1\vee S^1$ or even $\bigvee_i S^1$. And in the appendix B, we use prove the Nielsen-Schreier theorem which we will use this to study finitely index subgroups of $F_n$.
	\begin{theorem} For $ m $ and $ n $ non-negative integers, the rank of a free group is well-defined. That is, if $ F_n $, the free group of rank $ n $, is isomorphic to $ F_m $, the free group of rank $ m $, then $ m = n $.
		\vspace{-5 mm}
		\begin{proof}
			Suppose we have an isomorphism $\phi: F_m\to F_n$. Let's denote $\{a_1,\cdots, a_m\}$ be the free generating set of $F_m$ and $\{b_1,\cdots, b_n\}$ be the free generating set of $F_n$. As shown in the figure below, consider the canonical surjection $F_n\xrightarrow{Ab} \bigoplus_i \Z<b_i>$, we get a homomorphism $ \phi \circ Ab: F_m\to \bigoplus_i \Z<b_i>$. Then consider the canonical surjection $F_m\xrightarrow{Ab} \bigoplus_j \Z<a_j>$, because $ \bigoplus_i \Z<b_i> $ is abelian, there exists a unique homomorphism $F:\bigoplus_j \Z<a_j>\to \bigoplus_i \Z<b_i>$ s.t. the following diagram commutes. 
			\begin{displaymath}\xymatrix{
				F_m \ar[r]^{\phi}\ar[d]^{Ab}& F_n \ar[r]^{Ab}& \bigoplus_i \Z<b_i>\\
				\bigoplus_j \Z<a_j>\ar[rru]_F
			}\end{displaymath}
			Notice that because both $\phi$ and $Ab$ are surjective, so is the map $F$. So we get a surjective homomorphism $\Z^m \to \Z^n$. Now we use the map $F_n \xrightarrow{\phi^{-1}}F_m$, then by the same reason, we get a map $G:\bigoplus_i \Z<b_i>\bigoplus_j \Z<a_j>$ s.t. which is also a surjective homomorphism. So we have an isomorphism $\Z^m\to \Z^n$. So we have $m=n$. That is to say, the rank of a free group is well-defined.
		\end{proof}
	\end{theorem}
	\begin{theorem}
		For $ m $ and $ n $ non-negative integers,For which $ (m, n) $ does $ F_m $ have a subgroup isomorphic to $ F_n $?
		\vspace{-5 mm}
		\begin{proof}
			At first suppose $n=1$, then $F_1\backsimeq \Z$ is abelian. So we can let $m=1$ and then we have $F_1\backsimeq F_1$. Notice that if $m>1$, $F_m$ is not an abelian group, so it cannot have a subgroup isomorphic to $F_1$. \\
			Now let $n\geq 2$. Notice that for the map $F_1\to F_n$ sends the generating element of $F_1$ to the generating element of $F_n$. This is an injective homomorphism, so $F_1$ has a subgroup isomorphic to $F_n$. Similarly for the map $F_2\to F_n$ sends the two elements of the generating set of $F_2$ to the generating set of $F_n$. Notice that this is also an injective homomorphism, so $F_2$ has a subgroup isomorphic to $F_n$. Now we claim that $F_m$ has subgroup isomorphic to $F_2$ for all $m>1$. Let consider the fundamental group of $S^1\vee S^1$. Notice that it is just $F_2$. By \textbf{Theorem \ref{Injective}}, the map $p_*:\pi_1(\tilde{X},\tilde{x}_0)\to \pi_1(X,x_0)$ induced by the covering map is injective. So we suffice to show that the space $S^1\vee\cdots \vee S^1$ (with $m$ copies of $S^1$) is a covering space of $S^1\vee S^1$ (because the fundamental group of $S^1\vee\cdots \vee S^1$ is just $F_m$, this show $F_m$ has a subgroup isomorphic to $F_2$).\\
			For $m=2$, $ S^1\vee S^1$ is obviously a covering space of $ S^1\vee S^1$ (by the identity map as the covering map). \\
			For $m= 3$, consider the figure below,
			\begin{center}
				\includegraphics[scale = 0.4, angle = 0]{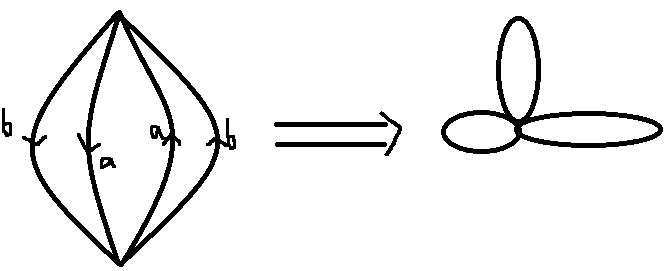}
			\end{center}
			This is a covering space of $S^1\vee S^1$ (by the covering map defined as the arrows). Notice that this is deformation retract to $S^1\vee S^1\vee S^1$ (when we quotient out one of edges as the figure above). So the fundamental group of $S^1\vee S^1\vee S^1$ ($F_3$) has a subgroup isomorphic $p_*$ to $F_2$.\\
			For $m>3$, similar as the above pattern, consider the figure below,
			\begin{center}
				\includegraphics[scale = 0.4, angle =0]{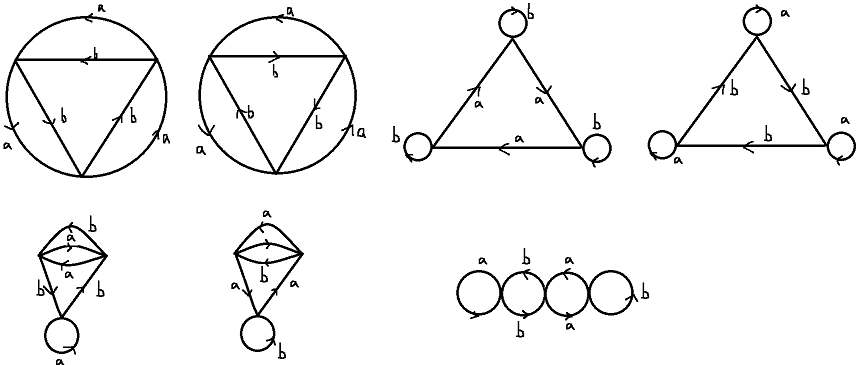}
			\end{center}
			This is a covering space of $S^1\vee S^1$ (by the covering map defined as the arrows). 
			Notice that this is deformation retract to $S^1\vee \cdots \vee S^1$ (as the figure above). So the fundamental group of $S^1\vee \cdots\vee S^1$ ($F_m$) has a subgroup isomorphic $p_*$ to $F_2$.\\
			Then because $F_2$ has a subgroup isomorphic to $F_n$, so $F_m$ has a subgroup isomorphic to $F_n$ for $m>1$.
			Then overall for $n\geq 2$, $F_m$ will directly have subgroup isomorphic to $F_n$ (for all $m>0$, when $m=0$, trivial).
		\end{proof}
	\end{theorem}
	\begin{theorem}
		For $ m $ and $ n $ non-negative integers, for which $ (m, n) $ does $ F_m $ have a normal subgroup isomorphic to $ F_n $?
		\vspace{-5 mm}
		\begin{proof}
			At first suppose $n=m=1$, that is trivial because they are isomorphism. If $n=1$, $ m>1 $, $F_m$ is not even subgroup isomorphic to a subgroup of $F_1$. If $m=1$ and $n>1$, then obviously $F_1$ is abelian and so the homomorphism sends the generating element of $F_1$ to a generating of $F_n$ sends $F_1$ an abelian subgroup isomorphic to $F_n$, so it is normal.\\
			Next, suppose $m,n>1$, we claim the \textbf{Theorem \ref{index}}: If $F$ is finitely generated free group and $N$ is a nontrivial normal subgroup of infinite index, then $N$ is not finitely generated. Now suppose $F_m$ has a subgroup isomorphic to the $F_n$, then the claim implies that $F_m$ can be isomorphism to a subgroup with finite index. Because a subgroup of a free group is free, the normal subgroup is also free and by part (1), the rank of that normal group is also $m$. And then by the Nielsen-Schreier formula \textbf{Theorem \ref{Nie}}, we have $m=1+e(n-1)$ where $e$ is the index of that normal subgroup. That is to say, a necessary condition that $F_m$ has a subgroup isomorphic to $F_n$ is $n-1 |m-1$. \\
			Next we show that if $n-1| m-1$, then $F_m$ has a normal subgroup isomorphic to $F_n$. If $n=m$, that is straightforward. Let $k(n-1)=m-1$ for $k\geq 2$. Let $S^1\vee\cdots \vee S^1$ with $m$ copies is the space $\tilde{X}$ and $S^1\vee\cdots \vee S^1$ with $m$ copies be the base space. As the figure below, we define a map $p:\tilde{X}\to X$ as the figure to make $\tilde{X}$ be a covering space of $X$.
			\begin{center}
				\includegraphics[scale = 0.4, angle = 0]{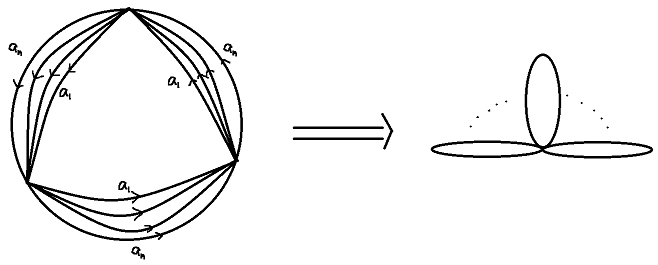}
			\end{center}
			Notice that for each $x\in X$ and each pair $\tilde{x},\tilde{x}'$ in the fiber of $x$, there is a deck transformation taking $\tilde{x}$ to $\tilde{x}'$ by the rotation map with feasible degree. Therefore, $\tilde{X}$ is a normal covering space. Then by \textbf{Theorem \ref{Normal}}, $p_*(\pi_1 (\tilde{X},\tilde{x}_0))$ is a normal subgroup of $\pi_1 (X,x_0)\backsimeq F_n$. That is to say $F_m$ has a normal subgroup isomorphic to $F_n$.
		\end{proof}
	\end{theorem}
	\begin{theorem}
		For $ m $ and $ n $ non-negative integers, for which $ (m, n) $ does $ F_m $ have a quotient group isomorphic to $ F_n $?
		\vspace{-5 mm}
		\begin{proof}
			Claim: $F_m$ has a quotient group isomorphic to $F_n$ iff $m\geq n$.\\
			To the forward direction, suppose we have a normal subgroup $N$ s.t. $F_m/N\backsimeq F_n$. Let's denote the isomorphism by $\phi: F_m/N\to F_n$ and denote the quotient map $\pi: F_m\to F_m/N$. So we have a (surjective) homomorphism $F_m \xrightarrow{\phi\circ \pi} F_n$. Similar as part (1), Let's denote $\{a_1,\cdots, a_m\}$ be the free generating set of $F_m$ and $\{b_1,\cdots, b_n\}$ be the free generating set of $F_n$. As shown in the figure below, denote $\psi:=\phi\circ \pi$ and consider the canonical surjection $F_n\xrightarrow{Ab} \bigoplus_i \Z<b_i>$, we get a homomorphism $ \psi \circ Ab: F_m\to \bigoplus_i \Z<b_i>$. Then consider the canonical surjection $F_m\xrightarrow{Ab} \bigoplus_j \Z<a_j>$, because $ \bigoplus_i \Z<b_i> $ is abelian, there exists a unique homomorphism $F:\bigoplus_j \Z<a_j>\to \bigoplus_i \Z<b_i>$ s.t. the following diagram commutes. 
			\begin{displaymath}\xymatrix{
				F_m \ar[r]^{\psi}\ar[d]^{Ab}& F_n \ar[r]^{Ab}& \bigoplus_i \Z<b_i>\\
				\bigoplus_j \Z<a_j>\ar[rru]_F
			}\end{displaymath}
			Notice that because both $\psi$ and $Ab$ are surjective, so is the map $F$. So we get a surjective homomorphism $\Z^m \to \Z^n$. That is to say, $m\geq n$.\\
			Conversely, suppose $m\geq n$.Then there is a natural surjection from $F_m\to F_n$ by sending the first $n$ generating element from $F_m$ to the generating set of $F_n$. Then by the 1st Isomorphism Theorem of groups, this surjection means $F_m$ have a quotient group isomorphic to $F_n$.
		\end{proof}
	\end{theorem}
	\begin{theorem}
		For $ m $ and $ n $ non-negative integers, list all the index 3 subgroups of $ F_2 $ and indicate which ones are conjugates. For each one, give a list of its free generators.
		\vspace{-5 mm}
		\begin{proof}
			Consider $X=S^1\vee S^1$. By the \textbf{Theorem \ref{sheet}}, the number of sheets of a path-connected covering space equals the index of $p_*(\pi_1(\tilde{X},x_0))$ in $\pi_1(X,x_0)$. The following is all the connected 3-sheeted covering spaces of $S^1\vee S^1$, up to isomorphism of covering spaces without basepoints. As shown below.
			\begin{center}
				\includegraphics[scale = 0.4, angle = 0]{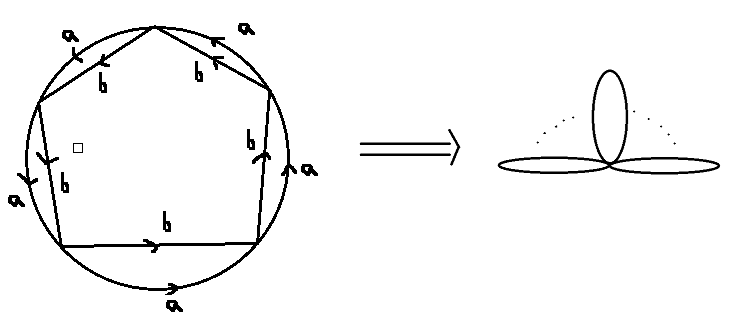}
			\end{center}
			By the Theorem 1.38 on page 67, there exists an one to one corresponding between the isomorphism classes of path-connected covering spaces and conjugacy classes. So none of $p_*(\pi_1(\tilde{X},\tilde{x}_0))$ for covering space $\tilde{X}$ above conjugates.
		\end{proof}
	\end{theorem}
	\begin{theorem}
		For $ m $ and $ n $ non-negative integers, does $ F_2 $ have any other rank four subgroups other than the ones listed above? If so, are any of them normal?
		\vspace{-5 mm}
		\begin{proof}
			Let's consider in two cases. At first suppose the subgroup has finite index, then by the Nielsen-Schreier formula (\textbf{Theorem \ref{Nie}}), the index of those subgroups has index 3. Then they are exactly $p_*(\pi_1(\tilde{X},\tilde{x}'))$ where $\tilde{X}$ are the covering spaces above. Notice that as the figure above, there are four of them are normal covering space (as the argument in part (3) for each $x\in X$ and each pair $\tilde{x},\tilde{x}'$ in the fiber of $x$, for that four covering space, there is a deck transformation taking $\tilde{x}$ to $\tilde{x}'$ by the rotation map with feasible degree. Therefore, $\tilde{X}$ is a normal covering space. so by the \textbf{Theorem \ref{Normal}}, $p_*(\pi_1(\tilde{X},\tilde{x}'))$ is a normal subgroup). Next consider the infinite index case. Consider the figure below. 			\begin{center}
				\includegraphics[scale = 0.4, angle = 0]{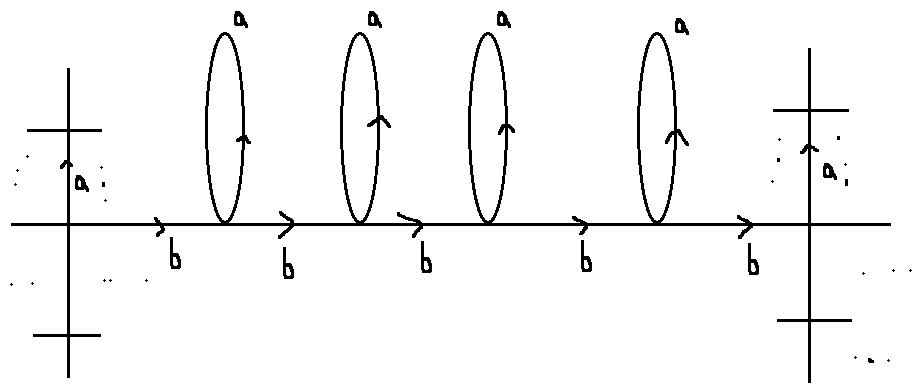}
			\end{center}
			It has rank four but it is infinitely sheeted, so it is infinite index. Notice that by \textbf{Theorem \ref{index}}, if $p_*(\pi_1(\tilde{X},\tilde{x}'))$ for that covering space above is a normal subgroup, then $p_*(\pi_1(\tilde{X},\tilde{x}'))$ cannot be finitely generated. That is to say, for all infinite index subgroup, it cannot be normal. 
		\end{proof}
	\end{theorem}
	\appendix

	\chapter{An introduction to Free groups and covering spaces}
	Let $S$ be a nonempty set, we denote $S^{-1}$ be any set disjoint from $S$ s.t. the cardinality are equal. Let's denote $s\in S$ its corresponding element $s^{-1}$ in $ S^{-1} $. Take a singleton set which is not contained in $S^\cup S^{-1}$ (let's denote it $\{e\}$). Then we define a word on $S$.
	\begin{definition}
		A word on $S$ is a sequence $(s_1,\cdots,s_n)$ where $s_i\in S\cup S^{-1}\cup \{e\}$ for $1\leq i\leq n$. A reduced word is a word s.t. 
		\vspace{-7mm}
		\begin{enumerate}
			\item $ s_{i+1}\neq s_i^{-1} $ for all $i$ with $s_i\neq e$
			\item If $s_k=e$ for some $k$, then the word is $e$.
		\end{enumerate}
		\vspace{-7mm}
	So then we define the set of reduced words as $F(S)$.
	\end{definition}
	\begin{theorem} $F(S)$ is a group under concatenation. We define $F(S)$ a free group.
	\end{theorem}
	\begin{theorem} (The Universal Property of free groups) 
		Let $G$ be a group, with a map $\psi: S\to G$, then there is a unique group homomorphism $\Psi:F(S)\to G$ s.t. the following diagram commutes:
		\begin{displaymath}\xymatrix{
			S \ar@{^{(}->}[r] \ar[dr]^{\psi}& F(S)\ar[d]^{\Psi}\\
			& G
		}\end{displaymath}	
	\end{theorem}
	\begin{corollary}
		$F(S)$ is unique up to isomorphism.
	\end{corollary}
	\begin{theorem}
	Every subgroup of a free group is free.
	\end{theorem}
	\begin{corollary}
		If $ F $ is a finitely generated free group and $ N $ is a nontrivial normal subgroup of infinite index, then $ N $ is not finitely generated.
	\end{corollary}
	Next we review some basic facts we need about covering spaces.
	\begin{definition}
		A \textbf{covering space} of a space $X$ is a space $ \tilde{X} $ together with a map $p:\tilde{X} \to X$ s.t. for each $x\in X$, there exists an open neighborhood $U$ in $X$ s.t. $p^{-1}(U)$ is a union of disjoint open sets in $\tilde{X}$, each of which is mapped homeomorphicaly onto $U$ by $p$.
	\end{definition}
	\begin{theorem} \label{Injective}The map $ p _* :\pi_1 ( \tilde{X}, \tilde{x}_0 ) \to  \pi_1 (X,x_0 ) $ induced by a covering space $ p:( \tilde{X}, \tilde{x}_0) \to  (X,x_0 ) $ is injective. The image subgroup $ p _* (\pi_1 ( \tilde{X}, \tilde{x}_0) )\in \pi_1 (X,x 0 ) $ consists of the homotopy classes of loops in $ X $ based at $ x_0 $ whose lifts to $ X $ starting at $ \tilde{x}_0 $ are loops.
	\end{theorem}
	\begin{theorem} \label{sheet} The number of sheets of a covering space $ p _* :\pi_1 ( \tilde{X}, \tilde{x}_0 ) \to  \pi_1 (X,x_0 ) $ with $ X $ and $ \tilde{X} $ path-connected equals the index of $ p _* (\pi_1 ( \tilde{X}, \tilde{x}_0) )\in \pi_1 (X,x 0 ) $.
	\end{theorem}
	\begin{theorem}\label{Hatcher about groups}
		If $ X $ is path-connected and locally path-connected, and semilocally
		simply-connected. Then there is a bijection between the set of basepoint-preserving
		isomorphism classes of path-connected covering spaces $ p : (\tilde{X},\tilde{x}_0) \to  (X,x_0) $ and the set of subgroups of $ \pi_1 (X,x_0 ) $, obtained by associating the subgroup $ p_* (\pi_1 ( \tilde{X},\tilde{x}_0 ) )$ to the covering space $ ( \tilde{X}, \tilde{x}_0 ) $. If basepoints are ignored, this correspondence gives a bijection between isomorphism classes of path-connected covering spaces $ p : \tilde{X} \to  X $ and conjugacy classes of subgroups of $\pi_1 (X,x_0)$. 
	\end{theorem}
	
	\begin{theorem}\label{Normal}
		Let $ p :( \tilde{X}, \tilde{x}_0 ) \to (X,x_0 ) $ be a path-connected covering space ofthe path-connected, locally path-connected space $ X $, and let $ H $ be the subgroup $ p_* (\pi 1 ( \tilde{X}, \tilde{x}_0 )) \subset \pi_1 (X,x_0 ) $. Then this covering space is normal iff $ H $ is a normal subgroup of $ \pi_1 (X,x_0 ) $.
	\end{theorem}
	
	The proofs of the above theorems can be found by \cite{Hatcher}, \cite{Dummit}. Overall, we will use that facts about free groups.
	
	\chapter{Nielsen-Schreier theorem}
	\begin{theorem}\label{index} If $F$ is finitely generated free group and $N$ is a nontrivial normal subgroup of infinite index, then $N$ is not finitely generated.
	\end{theorem}

	\begin{theorem}(Nielsen-Schreier theorem)\label{Nie}  If G is a free group on $ n $ generators, and $ H $ is a subgroup of finite index $ e $, then $ H $ is free of rank $ 1+e(n-1) $.
		\vspace{-5 mm}
		\begin{proof}
			Let's prove it using algebraic topology again. By the theorem above, we have an $e$-sheeted (path-connected) covering space (denote it $\tilde{X}$) of $X=\bigvee_i^n S^1$. 
			\begin{claim*} \normalfont 
				For $ X $ a finite CW complex and $ p :\tilde{X} \to X $ an $ e $-sheeted covering space, show that $ \chi (\tilde{X}) = e \chi (X) $.
				\vspace{-7 mm}
				\begin{proof}
					Consider the characteristic maps $\varphi_\alpha : D^k \to X$, because $D^k$ is contractible, we can lift it by the lifting criterion. Because $\tilde{X}$ is $n$-sheeted. There are exactly $e$ liftings of $\varphi_\alpha$ to $Y$. So for each $k$-cell $e^k$ in $X$, there exists $e$ $k$-cells in the lifted CW-structure of $\tilde{X}$. Now let $a_i$ be the number of $i$-cells of $\tilde{H}$, and $b_i$ be the number of $i$-cells in $X$. We have $a_i = e \cdot b_i$. And $ \chi(\tilde{H}) = \sum_{i = 0}^n (-1)^i a_i = \sum_{i = 0}^n (-1)^i e b_i = e \chi(X) $.
				\end{proof}
			\end{claim*}
		Let's figure out $\chi(X)$. Notice that the CW complex construction of $X$ is 1 0-cell and $n$ 1-cells. So $\chi(X)=1-n$. On the other hand, let's denote $\tilde{X}$ has $1$ 0-cell and $x$ 1-cells. Because $X$ is path-connected and orientable, so is it path-connected covering space $\tilde{X}$. So $H_1(\tilde{X})$ is free. Then by the Euler-Characteristic formula from \cite{Hatcher}, we have $\chi(\tilde{X})=\sum_{i = 0}^n (-1)^i \rank H_i (\tilde{X})=1-x$. That is to say $\chi(\tilde{X})=e\chi(X)=e(1-n)=1-x$. So $\rank H_1(\tilde{X})=x=1+e(n-1)$. I.e., $ H_1(\tilde{X})=\Z^{1+e(n-1)} $. Notice that $H_1$ is an abelianization of the fundamental group. We have $\pi_1(\tilde{X})\backsimeq \bigast_{1+e(n-1)}\Z=F_{1+e(n-1)}$.
		\end{proof}
	\end{theorem}
	\nocite{Bredon} \nocite{Dummit} \nocite{Hatcher} \nocite{Tu}
	\bibliographystyle{ieeetr}
	\bibliography{bib}
	
\end{document}